\documentclass[a4paper,12pt]{article}
\input psfig.sty

\usepackage{amsmath,amssymb}

\newtheorem {lemma}{Lemma}
\newtheorem {claim}{Claim}

\newtheorem {thm}{Theorem}
\newcommand{\rmd}{\,\mathrm{d}}

\def \a   {\alpha}
\def \b   {\beta}
\def \d   {\delta}
\def\ba {\bar\a}
\def\bb {\bar\b}
\def\bg {\bar\g}

\def \g   {\gamma}
\def \l   {{\lambda}}
\def \s   {{\sigma}}
\def \eps {\varepsilon}
\def\phi{\varphi}

\def\E{{\mathbb{E\,}}}
\def\Var{{\mathbb{V}{\sf ar}\,}}
\def\P{{\mathbb{P}}}
\def \R {{\mathbb{R}}}

\def\F{{\cal{F}}}
\def\limt{\lim_{t\to\infty}}
\def\limt0{\lim_{t\to 0}}

\def\|{\,|\,}
\def \eps {\varepsilon}
\newcommand{\BBox}{\rule{6pt}{6pt}}
\newcommand\Cox{$\hfill \BBox$ \vskip 5mm}
\def\bn#1\en{\begin{align*}#1\end{align*}}
\def\bnn#1\enn{\begin{align}#1\end{align}}
\def\proof{\noindent{\sf Proof.\ }}

\title {Random geometric subdivisions}

\author{Stanislav Volkov\footnote{Department of Mathematics,University of Bristol, BS8~1TW, U.K.\newline E-mail:~S.Volkov@bristol.ac.uk}
}
\begin {document}
\maketitle

\begin {abstract}
We study several models of random geometric subdivisions arising from the model of Diaconis and Miclo (2011).
In particular, we show that the limiting shape of an indefinite subdivision of a quadrilateral is a.s.\ a parallelogram. We also show that the geometric subdivisions of a triangle by angle bisectors converge (only weakly) to a non-atomic distribution, and that the geometric subdivisions of a triangle by choosing random points on its sides converges to a ``flat'' triangle, similarly to the result of Diaconis and Miclo (2011).
\end {abstract}


\section{Motivation}
The aim of this paper is to consider several models involving random subdivision of geometrical objects.
Markov chains involving geometry and polygons have been studied quite widely in the literature, see e.g.~\cite{Ding} and references therein; however, the main motivation of this paper comes from the paper by Diaconis and Miclo~\cite{DM}, who considered the following model, earlier also studied by Hough~\cite{HOU}. A non-degenerate triangle is divided by its three medians into $6$ smaller triangles, and one of these new triangles is chosen with equal probability, and (optionally) re-scaled, thus becoming the ``new'' triangle. This procedure is repeated indefinitely. It was shown that in some sense the limiting triangle will be ``flat'', that is the largest angle will converge to~$\pi$.

In the current paper we consider several generalizations of the above model.
In Section~\ref{quadri} we consider subdivisions of a quadrilateral by the lines connecting the middle points on the opposite sides; in Section~\ref{sec-bisec} we consider subdivision of a triangle using angle bisectors. Finally, in Section~\ref{sec-subtr} we consider a sequence of triangles obtained by randomly choosing a point on each of the sides and letting them be vertices of the new triangle.

\section{Random subdivision of quadrilateral}\label{quadri}
\setlength{\unitlength}{1mm}
\begin{picture}(160,80)(-30,0)
 \shortstack[c]{
\thicklines
\put(0, 0){\line(1, 0){80}}
\put(0, 0){\line(1, 6){10.75}}
\put(80, 0){\line(-1, 5){8}}
\put(72, 40){\line(-5, 2){61}}
\thinlines
\put(40, 0){\line(0, 1){53}}
\put(5.4, 32){\line(5, -1){71}}
 \put(1.6,1.6){$A$}
 \put(81,1.6){$B$}
 \put(73, 40){$C$}
 \put(14,64){$D$}
 \put(41,1.6){$E$}
 \put(77,19){$F$}
 \put(41,54){$G$}
 \put(0.8,32){$H$}
 \put(41,27){$M$}
}
\end{picture}
\\

We are given a convex quadrilateral $ABCD$. Let $E,F,G,H$ be the middle points of the sides $AB$, $BC$, $CD$, and $DA$ respectively. Let $M$ be the point of intersection of segments $EG$ and $FH$.
Now we replace $ABCD$ by one of the following four quadrilaterals $AEMH$, $EBFM$, $MFCG$ and $HMGD$
with equal probabilities. Suppose that we repeat this procedure indefinitely. What is the limiting {\it shape} of the quadrilateral obtained as the limit of this procedure?

\begin{thm}
Under the procedure described above, the limiting shape of the quadrilateral will be a parallelogram, and the rate of convergence is geometric.
(Note that the shape of parallelogram is ``invariant'' for the procedure).
\end{thm}

\proof Observe that $\vec{HF}=\frac 12(\vec{AB}+\vec{DC})$, and $\vec{HM}=\vec{MF}=\frac 12 \vec{HF}$ where $\vec{x}$ denotes the vector $x$.
Also suppose that when we replace the original quadrilateral by one of the four smaller ones, we rescale the smaller one twice thus making it bigger; this will not affect the shape. Let $\vec{u}_0=\vec{AB}$ and $\vec{v}_0=\vec{DC}$ be the vectors corresponding to the ``horizontal'' sides of $ABCD$, and $\vec{u}_1$, $\vec{v}_1$ be the corresponding vectors of the new quadrilateral obtained by subdivision. The crucial observation is that
 \bn
 \{\vec{u}_1,\vec{v}_1\}=\left\{\begin{array}{ll}
 \left\{\vec{u}_0,\frac 12\left(\vec{u}_0+\vec{v}_0\right)\right\} &\text{with probability }1/2;\\
 \left\{\vec{v}_0,\frac 12\left(\vec{u}_0+\vec{v}_0\right)\right\} &\text{with probability }1/2.
 \end{array}
 \right.
 \en
Similarly we can define $\left\{\vec{u}_n,\vec{v}_n\right\}$, $n=2,3,\dots.$

\setlength{\unitlength}{1mm}
\begin{picture}(160, 50)(-10,0)
  \thicklines
  \put(0, 0){\vector(1, 0){90}}
  \put(0, 0){\vector(1, 4){10}}
  \put(0, 0){\vector(2, 1){45}}
  \put(0, 0){\vector(1, 1){30}}
  \thinlines
  \put(10, 40){\line(2, -1){80}}
  \put(92,0){$U_0$}
  \put(12,42){$V_0$}
  \put(46,24){$X_1$}
  \put(31,31){$X_2$}
  \put(-4,-3){$O$}
\end{picture}
\\

Let us place all vectors $\vec{u}_n,\vec{v}_n$ at the origin $O$ and let $U_n$ and $V_n$ be the corresponding endpoints of these vectors. Let $X_n=U_n$ if $U_n\notin\{U_{n-1},V_{n-1}\}$ and $X_n=V_n$ if $V_n\notin\{U_{n-1},V_{n-1}\}$ (exactly one of these two statements must be true). Then we see that all points $X_n$ lie on the segment $U_0 V_0$. Moreover, $X_1$ lies exactly in the middle of $U_0 V_0$, $X_2$ with equal probabilities splits $U_0 X_1$ or $V_0 X_1$ in the middle, etc. As a result, we see that $X_n$ a.s.\ converges to a point $X_{\infty}$ which is uniformly distributed on the segment $U_0 V_0$. Taking into account that $|U_n V_n|=2^{-n} |U_0 V_0|$, we obtain a deterministic speed of convergence of $\vec{u}_n$ and $\vec{v}_n$ towards $\vec{u}_{\infty}:=\vec{OX_\infty}$:
 \bn
  \frac{|U_0 V_0|}{2^{n+1}} \le \max\{|\vec{u}_n-\vec{u}_{\infty}|,|\vec{v}_n-\vec{u}_{\infty}|\}\le \frac{|U_0 V_0|}{2^n}.
 \en
The analogous statement holds also for the ``vertical'' sides corresponding to $AD$ and $BC$, thus yielding the desired convergence towards a parallelogram.
 \Cox

\section{Random subdivision of triangle with angle bisectors}\label{sec-bisec}

\setlength{\unitlength}{1mm}
\begin{picture}(120,100)(-10,0)
\thicklines
\put(0, 0){\line(1, 0){120}}
\put(0, 0){\line(3, 5){60}}
\put(120, 0){\line(-3, 5){60}}
\thinlines
\put(0, 0){\line(2, 1){92}}
\put(120, 0){\line(-2, 1){92}}
\put(60, 100){\line(0, -1){100}}
\put(-3,2){$A$}
\put(121,2){$B$}
\put(62, 100){$C$}
\put(61,-4){$D$}
\put(95,45){$E$}
\put(22,45){$F$}
\put(61,37){$M$}
\put(5,1){\tiny{$\a/2$}}
\put(4,5){\tiny{$\a/2$}}
\put(110,1){\tiny{$\b/2$}}
\put(111,5){\tiny{$\b/2$}}
\put(55.5,90){\tiny{$\g/2$}}
\put(61,90){\tiny{$\g/2$}}
\put(50,2){\tiny{$\b+\g/2$}}
\put(61,2){\tiny{$\a+\g/2$}}
\put(30,46){\tiny{$\a+\b/2$}}
\put(26,41){\tiny{$\g+\b/2$}}
\put(81,46){\tiny{$\b+\a/2$}}
\put(86,41){\tiny{$\g+\a/2$}}
\end{picture}
\\

Unlike the median-subdivision model considered in~\cite{DM}, suppose that we subdivide the triangle $ABC$ by the three angle bisectors which intersect the sides $AB$, $BC$, $CA$ at points $E$, $F$, $G$ respectively, and let $M$ be the point of intersections of all angle bisectors, the centre of the inscribed circle. The replacement procedure now states that the triangle $ABC$ is replaced with probabilities $1/6$ by one of the following triangles: $ADM$, $DBM$, $BEM$, $ECM$, $CFM$, $FAM$. As before, the object of interest is the {\em shape} of the limiting triangle obtained by indefinite repetition of the above replacement procedure.

It turns out to be convenient to work with the angles of the triangle. If the original triangle has the angles $(\a,\b,\g)$, $\a+\b+\g=\pi$, then the new triangle will have one of the following six sets of angles:
 \bnn\label{eqabg}
 \left(\frac{\a}2,\g+\frac{\b}2,\frac{\a+\b}2\right), \ \left(\frac{\a}2,\b+\frac{\g}2,\frac{\a+\g}2\right),
 \nonumber \\
 \left(\frac{\b}2,\a+\frac{\g}2,\frac{\b+\g}2\right), \ \left(\frac{\b}2,\g+\frac{\a}2,\frac{\a+\b}2\right), \\
 \left(\frac{\g}2,\b+\frac{\a}2,\frac{\a+\g}2\right), \ \left(\frac{\g}2,\a+\frac{\b}2,\frac{\b+\g}2\right).\nonumber
 \enn
Obviously, we cannot expect convergence almost surely for this procedure (since e.g.\ any angle can be halved on the next step with probability $1/3$).

Observe that we can generate the sequence of triangles by {\it always} choosing the left-bottom triangle, that is by using the mapping $(\a,\b,\g)\to \left(\frac{\a}2,\g+\frac{\b}2,\frac{\a+\b}2\right)$, and then performing a random permutation of the set of three newly obtained angles. Formally, let $(\a_{n},\b_{n},\g_{n})$ denote the set of the angles of the $n$-th triangle, then
\bnn\label{eq_dyn}
(\a_{n+1},\b_{n+1},\g_{n+1})&=\s_n\left(\frac{\a_n}2,\g_n+\frac{\b_n}2,\frac{\a_n+\b_n}2\right)\nonumber\\
 &=\s_n\left(\frac{\a_n}2,\pi-\a_n-\frac{\b_n}2,\frac{\a_n+\b_n}2\right)
\enn
where $\s_n=\left\{\s^{(1)},\s^{(2)},\s^{(3)},\s^{(4)},\s^{(5)},\s^{(6)}\right\}$ is a random permutation of the set of three elements; $\s_n\big((a,b,c)\big)$ takes one of the six possible values
$$
\left\{ (a,b,c),\ (b,c,a),\ (c,a,b),\ (c,b,a),\ (b,a,c),\ (a,c,b)  \right\}
$$
with equal probability, and $\s_n$'s are i.i.d. The advantage of this method is that the distribution of the set of the angles of the $n$-th triangle is completely symmetric with respect to exchanges of its components, even though the components are not independent. Unlike the case of barycentric subdivision studied in~\cite{DM}, only the weak limit exists in our case.
\begin{thm}
The sequence $(\a_n,\b_n,\g_n)$ converges in distribution to some limit.
\end{thm}
\proof
Let $S=\{(\a,\b,\g): \ \a,\b,\g\ge 0,\ \a+\b+\g=\pi\}$ be the 2-simplex with the standard Euclidean distance. Let $f_i(u)$, $u\in S$, $i=1,\dots,6$ be the set of functions given by (\ref{eqabg}). Let $u=(\a,\b,\g)$ and $v=(\a+x,\b+y,\g+z)$, so that $x+y+z=0$.
Observe that all $f_i$ are Lipschitz and that
\bn
\sum_{i=1}^2 \log \frac{|f_i(u)-f_i(v)|}{|u-v|}
&= \log \sqrt{\frac {[x^2+y^2+[y+z]^2+xz][x^2+z^2+[y+z]^2+xy]}{4(x^2+y^2+z^2)^2} }
\\
&=\frac 12 \log \frac {(y^2+3z^2+3yz)(3y^2+z^2+3yz)}{16(y^2+z^2+yz)^2} \\
&\le\frac 12 \log \frac {(3y^2+3z^2+3yz)(3y^2+3z^2+3yz)}{16(y^2+z^2+yz)^2}=\log \frac 34 <0.
\en
The identical bound holds for $i=3,4$ and $i=5,6$. Therefore, if $Z_u$ denotes a set of the angles obtained from $u$ by random subdivision, we have
$$
 \E \left[\log \frac{|Z_u-Z_v|}{|u-v|}\right]
  =\sum_{i=1}^6 \frac 16 \log \frac{|f_i(u)-f_i(v)|}{|u-v|}\le \log \frac {\sqrt{3}}2 <0.
$$
Therefore, by Theorem~1 from~\cite{BE}, see also the proof of Lemma~5.1 in~\cite{DM}, the mapping $u\to Z_u$ is ergodic, that is there is a (unique) probability measure $\nu$ on $S$ such that for any starting configuration $u=(\a,\b,\g)\in S$ we have $Z^{(n)}_u\to \nu$ weakly. (Here $Z^{(n)}$ stands for the superposition of $n$ i.i.d.\ mappings $Z$.)
\Cox

Now let us try to get a handle on the distribution of the limiting triangle.
For the purpose of simplicity, and without loss of generality, assume from now on that $\a_n+\b_n+\g_n \equiv 1$ (as opposed to $\pi$).

%

We already know that $(\a_n,\b_n)$ converges to some pair $(\ba,\bb)$ in distribution. Since $\ba$ is bounded and $\bb$ and $\bg:=1-\ba-\bb$ have the same distribution as $\ba$, and $\ba\bb$, $\bb\bg$, $\bg\ba$ all have the same distribution as well (from the symmetry), we conclude
\bn
 \E[\ba+\bb+\bg]=1  &\Longrightarrow  \E \ba=\E \bb=\E \bg=\frac 13,\\
 \E[\ba+\bb+\bg]^2=1 &\Longrightarrow  3\E \ba^2+6 \E [\ba \bb]=1
\en
Also from~(\ref{eq_dyn}) we have
\bn
 \E \ba^2&=\frac 13\left[\frac 14 \,\E \ba^2+\E \left[\bar\g+\frac{\bb}2\right]^2+\frac14\, \E(\ba+\bb)^2\right].
\en
This in turn yields
\bn
 \E\ba^2=\frac 17,\ \E[\ba\bb]=\frac 2{21},\ \Var(\ba)=\frac 2{63},\ Cov(\ba,\bb)=-\frac 1{63}
\en
which sheds some light on the distribution of $\ba$ and the dependence between $\ba$ and $\bb$.

A more interesting and subtle statement about the joint distribution of $\ba,\bb,\bg$ is the following
\begin{thm}\label{th-bis}
Let $c_1,c_2,c_3$ be some real numbers, not all of which are $0$. Then distribution of the random variable $c_1 \ba+c_2 \bb+c_3 \bar\g$ does not have atoms.
\end{thm}

\begin{figure}[htb]
\centerline{\hbox{\psfig{figure=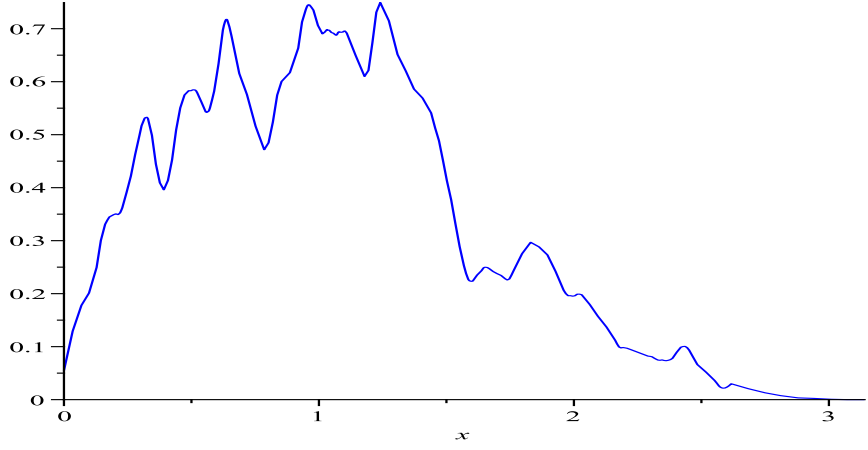,width=12cm,height=6cm}}}
\caption{Density of a randomly chosen angle.} \label{Fig1}
\end{figure}

In fact, we conjecture that the distribution of $\nu$ is continuous on the simplex, and so are the marginal distributions, e.g.\ the distribution of $\ba$.
Numerical simulations suggest that the pdf of a randomly chosen angle of a limiting triangle  looks like the one shown on Figure~\ref{Fig1}, which is obviously quite non-trivial.

\begin{figure}[htb]
\hbox{\psfig{figure=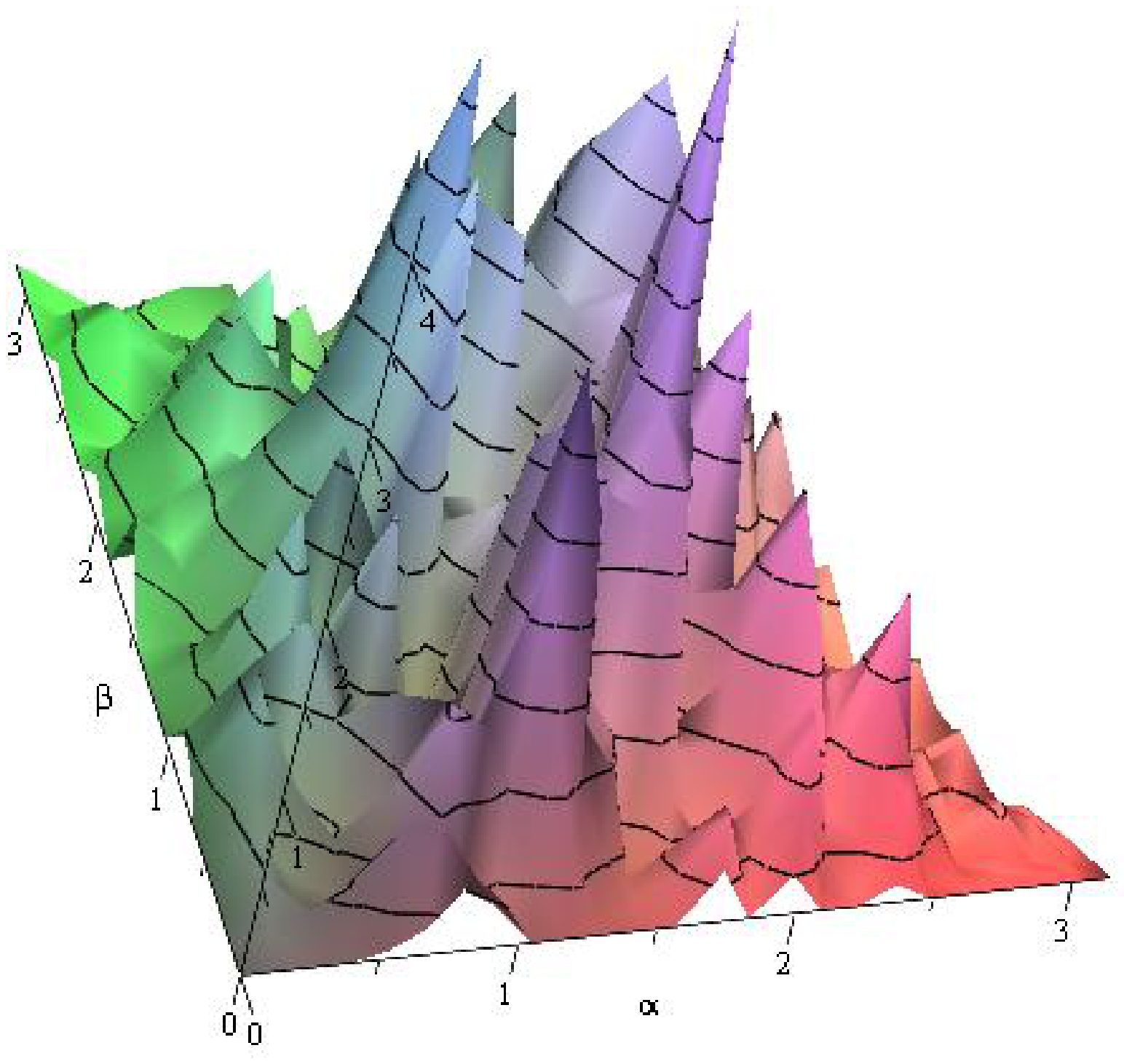,width=7cm,height=7cm}\psfig{figure=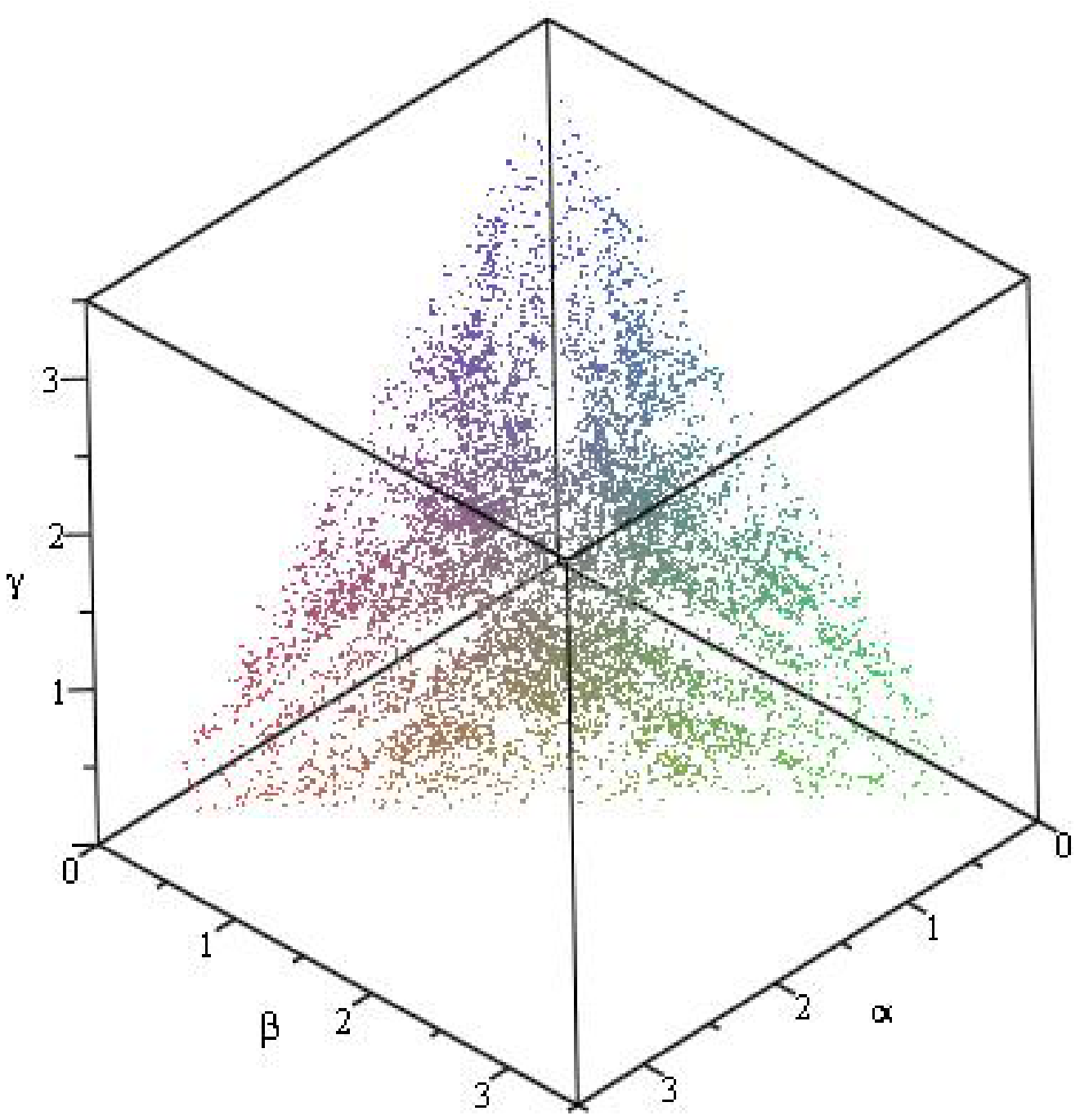,width=7cm,height=7cm}}
\caption{Sample density of $\nu$ on simplex $S$.} \label{Fig2}
\end{figure}

Due to the symmetry of the triple $(\ba,\bb,\bar\g)$ with respect to permutations and the fact that $\ba+\bb+\bar\g=1$, Theorem~\ref{th-bis} follows immediately from  the next statement.
\begin{lemma}
For any $c\in\R$ and $x\in\R$, $\P(\ba+c\bb=x)=0$.
\end{lemma}
\proof
Equation~(\ref{eq_dyn}) together with the weak convergence of $(\a_n,\b_n,\g_n)$ yield
\bnn\label{eq_next}
 \P(\ba+c\bb=x)&=\frac 16\left\{
 \P\left(\frac{\ba}2+c\left[1-\ba-\frac{\bb}2\right]=x\right)+
 \P\left(c\frac{\ba}2+\left[1-\ba-\frac{\bb}2\right]=x\right)\right.\nonumber\\
 &\ \ \ \ \ \ +
 \P\left(\frac{\ba}2+c\frac{\ba+\bb}2=x\right)+
 \P\left(c\frac{\ba}2+\frac{\ba+\bb}2=x\right)\nonumber\\
 &\ \ \ \ \ \ + \left. \P\left(\left[1-\ba-\frac{\bb}2\right]+c \frac{\ba+\bb}2=x\right)\right.\nonumber\\
 &\ \ \ \ \ \ +\left. \P\left(c\left[1-\ba-\frac{\bb}2\right]+\frac{\ba+\bb}2=x\right) \right\}.
\enn
In particular,
\bnn\label{eq_next1}
 \P(\ba=x)&=\frac 13\left\{\P(\ba/2=x)+\P(1-\ba-\bb/2=x)+\P([\ba+\bb]/2=x)\right\}
\enn
and if we set $x=1$, then
\bn
 \P(\ba=1)&=\frac 13\left\{\P(\ba/2=1)+\P(1-\ba-\bb/2=1)+\P([1-\bar\g]/2=1)\right\}\\
 &=\frac 13\left\{0+ \P(\ba=0,\bb=0)+0\right\}=\frac 13 \, \P(\bar\g=1)=\frac 13 \, \P(\ba=1),
\en
yielding
\bnn\label{eqba1}
 \P(\ba=1)=0.
\enn
Rewriting~(\ref{eq_next}) for $c\notin\{-1,1/2,2\}$ we have
\bnn\label{eq_nextt}
 \P(\ba+c\bb=x)&=\frac 16\left\{
 \P\left(\ba+\frac{c}{2c-1}\bb=?\right)+ \P\left(\ba+\frac{1}{2-c}\bb=?\right)\right.\nonumber\\
 &\ \ \ \ \ \ +\P\left(\ba+\frac{c}{1+c}\bb=?\right)+ \P\left(\ba+\frac{1}{1+c}\bb=?\right)\\
 &\ \ \ \ \ \  \left.+ \P\left(\ba+\frac{c-1}{c-2}\bb=?\right)+ \P\left(\ba+\frac{1-c}{1-2c}\bb=?\right) \right\}\nonumber
\enn
where the question marks stand for some non-random numbers. We can make sense of the expression above also for $c\in\{-1,1/2,2\}$  if we replace the fractions equal to infinity by zeros, due to the symmetry between $\ba$ and $\bb$.

For each $c$, let $x^*(c)$ be such that $\P(\ba+c\bb=x^*(c))=p^*(c)\equiv \max_{x\in\R}\P(\ba+c\bb=x)$; clearly such $x^*(c)$ must exist.
Suppose that the statement of the lemma does not hold; then $p^*=\sup_{c\in\R} p^*(c)>0$.
From~(\ref{eq_nextt}) and the symmetry between $\ba$ and $\bb$ it follows that
\bn
p^*(c)&=\P\left(\ba+c\bb=x^*(c)\right)\le\frac 16\left\{p^*+p^*+p^*+ p^*(c+1)+p^*+p^*\right\}\\
&=\frac 16\left\{5 p^*+ p^*(c+1)\right\},
\en
hence
\bn
p^*(c+1)=p^*(c)-6\eps \text{ whenever }p^*(c)> p^*-\eps.
\en
Fix a very small $\eps'>0$ and let $c_0$ be such that $p^*(c_0)> p^*-\eps'$.
Then for $c_i=c_0+i$, $i=1,2,\dots$, we have $p^*(c_i)>p^*-6^i\eps'$.

Let $A(c)=\{\omega:\ \ba+c\bb=x^*(c)\}$ and $Z_N=1_{A(c_0)}+1_{A(c_1)}+\dots+1_{A(c_{N-1})}$.
On one hand,
$$
\E Z_N=\sum_{i=0}^{N-1} \P(A(c_i))> Np^*- 0.2\cdot 6^{N} \eps'
$$
since $\P(A(c_i))=p^*(c_i)>p^*-6^i\eps'$.
On the other hand, for $1\le n \le N$,
$$
\E Z_N=\sum_{i=1}^{N} \P(Z_N\ge i)\le n-1+\P(Z_N\ge n),
$$
thus
$$
\P(Z_N\ge n)>N  p^* -(n-1) - 0.2\cdot 6^{N} \eps'.
$$
Suppose that
\bnn\label{eq_condN}
N  p^*> (n-1) + 0.2\cdot 6^{N} \eps'.
\enn
Then there is a subset of $A(c_0),\dots,A(c_{N-1})$ containing $n$ distinct elements, say $A(c_{i_1}),\dots,A(c_{i_n})$ for which the probability of $\cap_{m=1}^n A(c_{i_m})$ is strictly positive. We are going to make use of
\begin{claim}\label{clai1}
Suppose that $z_1,z_2,\dots,z_n$, $n\ge 3$, is a collection of $n$ distinct real numbers. Let $B=\bigcap_{i=1}^n A(z_i)$.
If $\P(B)>0$  then the sets $A(z_i)\setminus B$, $i=1,2,\dots,n$,  are disjoint and as a result
$$
\P\left(\bigcup_{i=1}^n A(z_i)\right)=\P(B)+\sum_{i=1}^n \P(A(z_i)\setminus B)=
 \sum_{i=1}^n \P(A(z_i))-(n-1) \P(B).
$$
\end{claim}
Proof of the claim. Since $\P(B)>0$, the system of equations
\bn\left\{\begin{array}{rc}
\ba+z_1\bb&=x^*(z_1)\\
\ba+z_2\bb&=x^*(z_2)\\
\dots\\
\ba+z_n\bb&=x^*(z_n)
\end{array}\right.
\en
must have a solution in $(\ba,\bb)$. This yields that for all distinct $i,j\in\{1,\dots,n\}$ we have
\bnn\label{eq_bb}
 \bb=\frac{x^*(z_i)-x^*(z_j)}{z_i-z_j}
\enn
and thus for any three distinct $i,j,k$
\bn
 \frac{x^*(z_i)-x^*(z_j)}{z_i-z_j}=\frac{x^*(z_i)-x^*(z_k)}{z_i-z_k}
\en
which in turn leads to $\omega\in A(z_i)\cap A(z_j) \Rightarrow \omega\in A(z_k)$ for all $k\Rightarrow \omega\in B$. Now the second statement of the claim is trivial.
\Cox

From Claim~\ref{clai1} it follows that for $B=\cap_{m=1}^n A(c_{i_m})$
$$
\P(B)\ge \frac{-1+\sum_{m=1}^n \P(A(c_{i_m}))}{n-1}>p^*-6^{N}\eps'-\frac 1{n-1}.
$$
Let
$
\eps''=6^{N}\eps'+\frac 1{n-1}.
$
From~(\ref{eq_bb}) we have that if
we set $y=\frac{x^*(c_{i_1})-x^*(c_{i_2})}{c_{i_1}-c_{i_2}}$
then $\P(\ba=y)=\P(\bb=y)\ge \P(B) \ge p^*-\eps''$.
At the same time, from~(\ref{eq_next1}) it follows
\bn
 p^*-\eps''&\le \P(\ba=y)=\frac 13\left\{\P\left(\ba=2y\right)+
 \P\left(\ba+{\bb}/2=2(1-y)\right) +
 \P\left(\ba+\bb=2y\right) \right\}
\\
&=\frac 13\left\{\P\left(\ba=2y\right)+
 \P\left(\ba+{\bb}/2=2(1-y)\right) +
 \P\left(\ba=1-2y\right) \right\}.
\en
Since each of the expressions on the RHS does not exceed $p^*$, we have $\min\{\P(\ba=2y),\P(\ba=1-2y)\}\ge p^*-3\eps''$. W.l.o.g.\ assume that  $2y\ge 1/2$.
Reiterating~(\ref{eq_next1}) and recalling~(\ref{eqba1}) (so that $\P(\ba\ge 1)=\P(\ba=1)=0$\,)  we obtain
\bn
 p^*-3\eps''&\le \P(\ba=2y)\\&=\frac 13\left\{\P\left(\ba=4y\right)+
 \P\left(\ba+\frac{\bb}2=2(1-2y)\right) +
 \P\left(\ba=1-4y\right) \right\}
\\
&\le \frac 13\left\{0+
 \P\left(\ba+{\bb}/2=2(1-2y)\right) +
 0 \right\}\le \frac {p^*}3
\en
leading to a contradiction, provided $\eps''<2p^*/9$. To finish the proof, we have to demonstrate that there exist $n,N,\eps'$ instantaneously satisfying  this condition as well as~(\ref{eq_condN}).
Indeed, fix an $n$  so large that $n-1>9/p^*$. Next, let $N$ be  larger than $(n-1)/p^*$.
Finally, set
$$
\eps'=\frac{\min\left\{Np^*-(n-1),p^*/9\right\}}{6^{N}}>0.
$$
It is easy to see that~(\ref{eq_condN}) is fulfilled and moreover
$$
\eps''=6^N \eps'+\frac 1{n-1}<\frac{p^*}9+\frac{p^*}9=\frac{2p^*}9.
$$
\Cox

\section{Random subtriangle of triangle }\label{sec-subtr}

\setlength{\unitlength}{0.5mm}
\begin{picture}(120,105)(-20,0)
\thicklines
\put(0, 0){\line(1, 0){120}}
\put(0, 0){\line(3, 5){60}}
\put(120, 0){\line(-3, 5){60}}
\thinlines
\put(40, 0){\line(2, 1){61}}
\put(101, 31){\line(-4, 1){71}}
\put(40, 0){\line(-1, 5){9.8}}
\put(-3,2){$A$}
\put(121,2){$B$}
\put(62, 100){$C$}
\put(40,-7){$C_1$}
\put(105,30){$A_1$}
\put(19,48){$B_1$}
\end{picture}
\\ \\
Now on each side of the triangle we randomly (independently and uniformly) choose points $A_1$, $B_1$, $C_1$.
The new triangle is now formed by these three points ($A_1 B_1 C_1$). Repeat this procedure indefinitely. What is the limit of the shape?

Let $a_0=|BC|$, $b_0=|CA|$, and $c_0=|AB|$. Also let $\xi_a=|B A_1|/|BC|$, $\xi_b=|C B_1|/|CA|$,  and $\xi_c=|A C_1|/|AB|$.
According to our assumption, $\xi_i$, $i=a,b,c$, are independent uniform $[0,1]$ random variables.
Also let $a_1=|B_1 C_1|$, $b_1=|C_1 A_1|$, and $c_1=|A_1 B_1|$ be the sides of the new obtained triangle.

\begin{thm}\label{thrtr}
The limiting triangle shape is a.s.\ flat (i.e., maximum angle converges to $\pi$). Moreover,
for any $c< 1+\frac{\log 4}3-\frac{\pi^2}9=0.365\dots$
$$
 y_n\le e^{-c n} \mbox{ for all sufficiently large $n$}
$$
where $y_n$ is the ratio of the height of the triangle corresponding to its largest side and the length of this side.
\end{thm}

Similarly to~\cite{DM}, let us rescale the triangle such that the largest side's length (say, $AB$) is $1$, and fit this side on the coordinate plane such that $A=(0,0)$, $B=(1,0)$. Also suppose that $C$ lies in the upper plane and $|AC|\ge |BC|$. Let $C=(x,y)$, then $x\in [1/2,1]$, $y\ge 0$, and the pair $(x,y)$ completely characterizes the shape of the triangle $ABC$. We also have
\bn
A_1=\left(1-(1-x)\xi_a,y\xi_a\right),\ B_1=\left(x(1-\xi_b),y(1-\xi_b)\right),\ C_1=\left(\xi_c,0\right).
\en
The new side lengths will then satisfy
\bn
a_1^2&=\left[x(1-\xi_b)-\xi_c\right]^2+\left[y(1-\xi_b)\right]^2\\
b_1^2&=\left[(1-(1-x)\xi_a)-\xi_c\right]^2+\left[y\xi_a\right]^2
\\
c_1^2&=\left[(1-(1-x)\xi_a)-x(1-\xi_b)\right]^2+\left[y(1-\xi_a-\xi_b)\right]^2
\en
and according to the standard formulas, the area of the triangle with vertex coordinates $A_1$,  $B_1$, and  $C_1$ equals
$$
\Delta:=\frac y2 \left[ \xi_a  \xi_b \xi_c +(1-\xi_a)(1-\xi_b)(1-\xi_c) \right]
$$
Thus the new value of $y$ is now
$$
y_1=\frac{2\Delta}{\max\{a_1,b_1,c_1\}^2}=y\cdot \frac{\xi_a  \xi_b \xi_c +(1-\xi_a)(1-\xi_b)(1-\xi_c)}{\max\{a_1^2,b_1^2,c_1^2\}}
$$
and also
\bn
x_1=\left\{\begin{array}{ll}
\frac{a_1^2+b_1^2-c_1^2}{2\max\{a_1,b_1\}^2} \text{ if } c_1=\min\{a_1,b_1,c_1\};\\
\frac{a_1^2+c_1^2-b_1^2}{2\max\{a_1,c_1\}^2} \text{ if } b_1=\min\{a_1,b_1,c_1\};\\
\frac{b_1^2+c_1^2-a_1^2}{2\max\{b_1,c_1\}^2} \text{ if } a_1=\min\{a_1,b_1,c_1\}.
\end{array}\right.
\en
which can be summarized as
\bn
y_1&=y\cdot r,\\
x_1&=\frac{a_1^2+b_1^2+c_1^2-2\min\{a_1^2,b_1^2,c_1^2\}}{2\max\{a_1^2,b_1^2,c_1^2\}}
\en
 where
\bn
r&=r(x,y)= \frac{R}{\max\{a_1^2,b_1^2,c_1^2\}}
\en
and
\bn
R=R(\xi_a,\xi_b,\xi_c)= \xi_a \xi_b \xi_c +(1-\xi_a)(1-\xi_b)(1-\xi_c).
\en

Observe that the denominator of $r(x,y)$ is increasing in $y$; hence $r(x,y)\le r(x,0)$.

From now on let us assume $y=0$. For $y=0$, we have that points $A_1$, $B_1$, $C_1$ all lie on the horizontal axes with coordinates
$\mu=1-(1-x)\xi_a$, $\nu=x(1-\xi_b)$ and $\xi_c$ respectively; we always have $0\le \nu \le \mu \le 1$
since $\nu\in[0,x]$ and $\mu\in [x,1]$. Consequently,
\bn
 S(x;\xi_a,\xi_b,\xi_c)&:=\max\{a_1,b_1,c_1\}\\
 &=\max\{\mu-\xi_c,\mu-\nu,\xi_c-\nu\}=\left\{\begin{array}{ll}
 \mu-\xi_c &\mbox{ if } \xi_c<\nu;\\
 \mu-\nu &\mbox{ if } \nu\le \xi_c\le \mu;\\
 \xi_c-\nu &\mbox{ if } \xi_c>\mu.
 \end{array}\right.
\en
Therefore,
\bn
 \E [r(x,0)\| \xi_a,\xi_b]&=\int_0^1 \frac{R(\xi_a,\xi_b,\xi_c)}{[\max\{a_1,b_1,c_1\}]^2} \rmd \xi_c
=I_1+I_2+I_3\\
 &=\int_0^{\nu} \frac{R(\xi_a,\xi_b,\xi_c)}{(\mu-\xi_c)^2} \rmd \xi_c
 +\int_{\nu}^{\mu} \frac{R(\xi_a,\xi_b,\xi_c)}{(\mu-\nu)^2} \rmd \xi_c
 +\int_{\mu}^1 \frac{R(\xi_a,\xi_b,\xi_c)}{(\xi_c-\nu)^2} \rmd \xi_c.
\en
Easy algebra gives
\bn
I_1&=
(\xi_a+\xi_b-1) \ln  \left(\frac{ \xi_b x+(1-x)(1-\xi_a)}  {1-\xi_a+\xi_a x}\right)
  +\frac { (1-\xi_b)\xi_a x} {1-\xi_a+\xi_a x},  \\
I_2&= \frac{\xi_a+1-\xi_b}2,\\
I_3&=(\xi_a+\xi_b-1) \ln  \left(\frac  {1-x+\xi_b x} {\xi_b x+(1-x)(1-\xi_a)}  \right)
+\frac{\xi_a(1-\xi_b)(1-x)}{1-x+\xi_b x},
\en
hence
\bn
 \E [r(x,0)\| \xi_a,\xi_b]&=(\xi_a+\xi_b-1) \ln  \left(\frac  {1-x+\xi_b x}{1-\xi_a+\xi_a x}  \right)
 + \frac{\xi_a+1-\xi_b}2
 \\&
 -\frac{\xi_a (1-\xi_b)(x^2(\xi_a+\xi_b-1)+ \xi_a(1-2x)-1)}{ (1-\xi_a+\xi_a x)(1-x+\xi_b x) }.
\en
As a result,
\bn
 \E [r(x,0)\| \xi_a]&=\int_0^1 \E [r(x,0)\| \xi_a,\xi_b] \rmd \xi_b
  \\
  &=\frac 1{2x^2(1-\xi_a+\xi_a x)} \left[ x(x+1-\xi_a(1-x)(2x+3-\xi_a(2-x)))\right.
  \\
&+\left.(1-2\xi_a)(1-\xi_a+\xi_a x)((1-x^2)\log(1-x) +x^2\log(1-\xi_a+\xi_a x)) \right]
\en
and thus
\bn
 \E [r(x,y)]\le \E [r(x,0)]&=\int_0^1 \E [r(x,0)\| \xi_a] \rmd \xi_a  =1.
\en

Let $(x_0,y_0)$, $(x_1,y_1)$, $(x_2,y_2)$, $\dots$ be the sequence of coordinates corresponding to the sequence of subdivided triangles. Assume $y_0>0$, then with probability $1$ we have $y_n>0$ for all $n$. Let $\F_n$ be the sigma-field generated by $\{(x_i,y_i),i=0,1,\dots,n\}$.
We have just established that
$$
\E[y_{n+1} \|\F_n]=y_n \E[r(x_n,y_n)]\le y_n.
$$
Thus $y_n\in [0,1]$ is a supermartingale which must converge a.s.\ to some limit $y_{\infty}$.
Let $E=\{\xi_a<0.1,\ \xi_b>0.9\}$, then $\P(E)=0.01>0$.
On the event $E$ we have $c_1\ge 0.8$, hence
$$
r(x,y)= \frac{\xi_a \cdot [\xi_b \xi_c]+(1-\xi_b)\cdot[(1-\xi_a)(1-\xi_c)]}{\max\{a_1,b_1,c_1\}^2}
\le \frac{2\cdot 0.1}{c_1^2}<\frac {1}{3}.
$$
Consequently, $\P(y_{n+1}<y_n/3\| \F_n)>0.01$ which implies $y_{\infty}=0$ a.s.

Note that
\bn
 \E \left[\log R\left(\xi_a,\xi_b,\xi_c\right) \right]=\int_0^1 \int_0^1  \int_0^1 \log (R) \rmd\xi_c\rmd\xi_b\rmd\xi_a
 =\frac{\pi^2}9-\frac 83=-1.57\dots
\en
and
\bnn\label{eq_kappa}
 \E [\log S(x;\xi_a,\xi_b,\xi_c)]&=
 \int_0^1 \rmd \xi_a \int_0^1 \rmd \xi_b \left[
 \int_0^{\nu} \log(\mu-\xi_c) \rmd \xi_c
 +\int_{\nu}^{\mu} \log (\mu-\nu) \rmd \xi_c\right.\nonumber\\
  &\left.+\int_{\mu}^1 \log(\xi_c-\nu) \rmd \xi_c\right]=
 -\frac{5}{6}-\frac {x^3\log(x)+(1-x)^3\log(x)}{3x(1-x)}\nonumber\\
  &
 \ge \frac{\log 4 -5}{6}=-0.602\dots=:-\kappa
\enn

Observe that
\bn
\log y_{n} -\log y_{n-1}&=\log r(x_{n-1},y_{n-1})\le \log r(x_{n},0)\\
&=\log R\left(\xi_a^{(n)},\xi_b^{(n)},\xi_c^{(n)}\right)
 -2\log S\left(x_{n-1};\xi_a^{(n)},\xi_b^{(n)},\xi_c^{(n)}\right)\\
&= \rho_n + 2 \s_n
\en
where $\rho_n=\log R\left(\xi_a^{(n)},\xi_b^{(n)},\xi_c^{(n)}\right)$ are i.i.d.\ with mean $\frac{\pi^2}9-\frac 83$ and $\s_n=-\log S\left(x_{n-1};\xi_a^{(n)},\xi_b^{(n)},\xi_c^{(n)}\right)$ are some $\F_n-$adapted random variables respectively.
Since $0\le S(\cdot) \le 1$, we have $\s_n\ge 0$, and also $\E (\s_n\|\F_{n-1})\le \kappa$ due to~(\ref{eq_kappa}). Additionally, we have for any positive $z$
\bn
\P(\s_n &\ge z\| \F_{n-1})=\P\left(S\left(x_{n-1};\xi_a^{(n)},\xi_b^{(n)},\xi_c^{(n)}\right) \le e^{-z}\| \F_{n-1}\right)
 \\
 &\le \P(\mu-\nu \le e^{-z} \| \F_{n-1})=\P(x_{n-1} \xi_b +(1-x_{n-1})(1-\xi_a) \le e^{-z} \| \F_{n-1})\\
 &\le\P(1_{\{x_{n-1}>1/2\}} x_{n-1} \xi_b +1_{\{x_{n-1}\le 1/2\}}(1-x_{n-1})(1-\xi_a)] \le e^{-z} \| \F_{n-1})\\
 &\le \P\left( \frac{\xi_b}2 \le e^{-z} \|\F_{n-1}\right)=2 e^{-z}
\en
since $1-\xi_a$ has the same uniform $(0,1)$ distribution as $\xi_b$.

The following statement is related to exponential inequalities involving martingales; however since we did not find it in the form we needed, we present its short proof later.
\begin{lemma}\label{lem_martineq}
Let $\F_n$ be an increasing family of $\sigma-$fields, and  $\s_n$ be $\F_n$-adapted random variables, possibly unbounded. Suppose that
\bn
 \E[\s_n\| \F_{n-1}]\le \kappa
\en
and for some $\bar z\ge 0$ and $c>0$ we have
\bn
 \P(|\s_n|\ge z\| \F_{n-1}]\le e^{-c z} \mbox{ for all }z\ge \bar z.
\en
Then
\bn
 \limsup_{n\to\infty} \frac{\sum_{i=1}^n \s_i}{n} \le \kappa \mbox{ a.s.}
\en
\end{lemma}

Now from Lemma~\ref{lem_martineq} and the strong law of large numbers it follows that
\bn
 \limsup_{n\to\infty} \frac{\log y_n}n&\le
  \limsup_{n\to\infty} \frac{\sum_{i=1}^n \rho_n}n
  +
  2 \limsup_{n\to\infty} \frac{\sum_{i=1}^n \s_n}n\\
  &\le  \E \rho_n+2\kappa=\frac{\pi^2}9-1-\frac{\log 4}3=-0.365\dots
\en
This established the result of Theorem~\ref{thrtr}.

\begin{thm}\label{thx}
$x_n$ converges to Uniform $[1/2,1]$ distribution.
\end{thm}
Proof.
Since $y_n\to 0$ a.s., and $x_{n+1}$ is continuous in $y_n$ near $0$, it will suffice to show that, given $y_n=0$, $x_{n+1}$ has $U[1/2,1]$ distribution. This will follow, in turn, from the following statement: put $\mu,\nu, \xi_c$ in an increasing order, and denote the resulting values $0 \le x^{(1)} <x^{(2)}<x^{(3)}\le 1$, then $\chi:=\frac{x^{(2)}-x^{(1)}}{x^{(3)}-x^{(1)}}$ has $U[0,1]$ distribution.
Indeed, for any $z\in(0,1)$ we have
\bn
\P(\chi\le z)&=\P(\chi\le z,\ \xi_c<\nu)+\P(\chi\le z,\ \nu\le \xi_c\le\mu)+\P(\chi\le z,\ \xi_c>\mu)\\
&=(I)+(II)+(III).
\en
We have
\bn
 (I)&=\int_0^1 \rmd\xi_b \int_0^1 \rmd\xi_a \int_0^{\nu}1_{\left\{ \frac {\xi_{a}x-\xi_{c}}{1-\xi_{b}(1-x)-\xi_{c}}<z
 \right\}} \rmd \xi_c
 =\left\{\begin{array}{ll}
 \frac { ( 3z - xz^2 - zx -x ) x}{6z (1-x) } &\text{ if }x<z;\\
 \frac { ( 3x - zx^2 - xz -z ) z}{6x (1-z) } &\text{ if }x\ge z,
 \end{array}\right.
\\ \\
 (II)&=\int_0^1 \rmd\xi_a \int_0^1 \rmd\xi_b \int_{\nu}^{\mu}1_{\left\{ \frac {\xi_c-\xi_a x}{1-\xi_b(1-x) -\xi_a x}<z \right\}} \rmd \xi_c=z/2,
\\ \\
 (III)&=\int_0^1 \rmd\xi_a \int_0^1 \rmd\xi_b \int_{\mu}^{1}1_{\left\{
  \frac {1-\xi_{b}(1-x)-\xi_{a}x}{\xi_{c}-\xi_{a}x}<z \right\}} \rmd \xi_c
 \\
 &=\left\{\begin{array}{ll}
   \frac {3z^2+z^2 x^2-3z^2 x + z x^2-3 z x+x^2}{6z(1-x) }&\text{ if }x<z;\\
   \frac { (1-x)^2 z^2}{6(1-z)x}
 &\text{ if }x\ge z.
 \end{array}\right.
\en
As a result, $\P(\chi\le z)=(I)+(II)+(III)=z$ and thus $\chi$ has a uniform $[0,1]$ distribution.
\Cox

\subsection{Proof of Lemma~\protect{\ref{lem_martineq}}}
Proof.
The result follows from the Borel-Cantelli lemma, once we  establish that for any $\eps>0$ there is a $\d=\d(\eps)>0$ such that
\bnn\label{eq4BC}
\P\left(\sum_{i=1}^n\s_i\ge (\kappa+\eps)n\right)\le e^{-\d n}.
\enn
First, by Markov inequality, for any $\l>0$
\bn
\P\left(\sum_{i=1}^n \s_i \ge (\kappa+\eps)n\right)
&=\P\left(e^{\l \sum_{i=1}^n \s_i} \ge e^{\l n(\kappa+ \eps)} \right) \le
e^{-\l n(\kappa +\eps)} \E\left[e^{\l \sum_{i=1}^n \s_i}\right]\\
&= e^{-\l n(\kappa +\eps)} \E \E\left[e^{\l \sum_{i=1}^n \s_i}\| \F_{n-1}\right]\\
&=e^{-\l n(\kappa +\eps)} \E \left[e^{\l \left(\sum_{i=1}^{n-1} \s_i\right)} \E[e^{\l \s_n}\| \F_{n-1}]\right].
\en
Now to show~(\ref{eq4BC}) by induction, it suffices to demonstrate that there are $\l>0$ and $\d>0$ such that $\E[e^{\l \s_n}\| \F_{n-1}]<e^{\l(\kappa+\eps)-\d}$ for all $n$.

Indeed, (see e.g.~\cite{DUR}, Lemma 5.7 in Section~1)
\bn
 \E\left[\frac{\l^k |\s_n|^k}{k!}\| \F_{n-1}\right]
 &=\frac{\l^k}{k!} \int_{0}^{\infty} kz^{k-1} \P(|\s_n|>z\| \F_{n-1}) \rmd z\\
 &\le \frac{\l^k}{k!} \left[ {\bar z}^k +\int_{\bar z}^{\infty} kz^{k-1} e^{-cz} \rmd z\right]\\
 &\le \frac{\l^k}{k!} \left[ {\bar z}^k +k! c^{-k}\right]= \frac{(\l \bar z)^k}{k!}+\left(\frac{\l}{c}\right)^k,
\en
therefore, assuming $\l<\min(c,1)$, by Dominated convergence theorem we have
\bnn\label{eq_elam}
 \E [e^{\l \s_n}\| \F_{n-1}]&=1+\l\E[\s_n\|\F_{n-1}]+\sum_{k=2}^{\infty} \frac{\l^k \E[\s_n^k\|\F_{n-1}]}{k!}
  \nonumber\\
 &\le 1+\l\kappa+\sum_{k=2}^{\infty} \E\left[\frac{\l^k |\s_n|^k}{k!}\|\F_{n-1}\right]\\
 &\le 1+\l\kappa+\sum_{k=2}^{\infty} \left(\frac{\l}{c}\right)^k+ \sum_{k=2}^{\infty}\frac{(\l \bar z)^k}{k!}\nonumber
 \\
 &\le 1+\l\kappa+\frac {\l^2}{c(c-\l )}+\l^2[e^{\bar z}-1-\bar z]\nonumber
\enn
Now by choosing $\l>0$ sufficiently small, we can make~(\ref{eq_elam}) smaller than $e^{\l(\kappa+\eps)-\d}>1+\l \kappa +(\l\eps-\d)$
for some $\d>0$, thus finishing the proof.
\Cox


\end{document}